\documentstyle[a4, pslatex, mathrsfs]{amsart}

\renewcommand{\cal}{\mathscr}
\newcommand{\catqot}{/\hskip-3pt/}
\newcommand{\C}{{\Bbb C}}

\newcommand{\E}{{\cal E}}

\newcommand{\GL}{\mathop{\rm GL}}

\newcommand{\Hom}{\mathop{\rm Hom}}

\renewcommand{\Im}{\mathop{\rm Im}}

\renewcommand{\L}{{\cal L}}

\renewcommand{\O}{{\cal O}}
\renewcommand{\P}{{\Bbb P}}

\newcommand{\SL}{\mathop{\rm SL}}
\newcommand{\Spec}{\mathop{\rm Spec}}

\renewcommand{\tilde}{\widetilde}

\newcommand{\la}{\lambda}
\newcommand{\lra}{\longrightarrow}

\newcommand{\p}{\prime}
\newcommand{\q}{\quad}

\renewcommand{\phi}{\varphi}
\newcommand{\rk}{\mathop{\rm rk}}

\renewcommand{\theta}{\vartheta}

\theoremstyle{plain}
\newtheorem{Thm}{Theorem}[subsection]
\newtheorem{Cor}[Thm]{Corollary}
\newtheorem{Prop}[Thm]{Proposition}
\newtheorem{Lem}[Thm]{Lemma}
\newtheorem*{MThm}{Theorem}

\theoremstyle{remark}
\newtheorem{Rem}[Thm]{Remark}

\begin{document}

\pagestyle{myheadings}
\title{The equivalence of Hilbert and Mumford stability for vector bundles}
\author{Alexander Schmitt$^*$}
\thanks{$^*$ Supported by a grant of the Emmy Noether Institute at Bar-Ilan
University.}
\address{Bar-Ilan University, Department of Mathematics \& Computer
Science, Ramat Gan, 52900, Israel}
\address{Universit\"at GH Essen, FB6 Mathematik \& Informatik, D-45117 Essen,
Deutschland\newline}

\email{schmitt1@@cs.biu.ac.il}


\begin{abstract}
We prove  
the equivalence of the notions of Hilbert (semi)stability
and Mumford (semi)stability for vector bundles on smooth curves for
arbitrary rank.
\end{abstract}

\maketitle
\markboth{Alexander Schmitt}{HILBERT STABILITY}

\section*{Introduction}
Consider the following setup. Let $d$, $g$, and $r$ be fixed positive 
integers, $W$ a complex vector space of dimension $p:=d+r(1-g)$,
and ${\frak G}=G(W,r)$ the Grassmannian of $r$-dimensional quotients
of $W$. On ${\frak G}$, there is the universal quotient
$$
W_{\frak G}\otimes\O_{\frak G}\lra E_{\frak G}
$$
which induces a surjection $\bigwedge^r W\lra \bigwedge^r E_{\frak G}$,
defining the Pluecker embedding 
${\frak G}\hookrightarrow \P(\bigwedge^rW)$. Now, let $C\hookrightarrow
{\frak G}$ be a smooth curve of genus $g$ such that $E_C:=E_{{\frak G}|C}$
has degree $d$. Then, $C$ gets embedded into $\P(\bigwedge^rW)$ as
a curve with Hilbert polynomial $P(m)=\chi((\bigwedge^rE_C)^{\otimes m})
=dm+(1-g)$. From the restriction $W\otimes\O_C\lra E_C$ of the universal
quotient to $C$ we derive homomorphisms $\bigwedge^r W\otimes\O_C\lra
\bigwedge^r E_C$, and for all
$m\ge 1$
$$
\psi^m_C\colon S^m\bigwedge^rW\lra 
H^0\bigl((\bigwedge^r E_C)^{\otimes m}\bigr).
$$
If $\psi^m_C$ is surjective (as will be the case for large $m$) and
$h^0((\bigwedge^r E_C)^{\otimes m})\allowbreak =P(m)$, this yields
$$
\phi_C^m:=\bigwedge^{P(m)} \psi_C^m\colon \bigwedge^{P(m)} 
(S^m\bigwedge^rW)\lra \C.
$$
We call $C$ or, abusively, $E_C$ \it $m$-Hilbert 
(semi/poly)stable\rm, if $\psi_C^m$ is surjective, 
$h^0((\bigwedge^r E_C)^{\otimes m})=P(m)$, and the point $\phi_C^m$
in $\P(\bigwedge^{P(m)} 
(S^m\bigwedge^rW))$ is (semi/poly)stable w.r.t.\ the natural action
of $\SL(W)$ on that space, and \it Hilbert (semi/poly)stable\rm ,
if it is $m$-Hilbert (semi/poly)stable for all $m$ sufficiently large.
This is now a new stability concept for the vector bundle
$E_C$ entering in competition to classical Mumford stability.
It goes back to Gieseker and Morrison 
(\cite{G1}, \cite{G2}). Its main motivation
is to obtain an alternative compactification, called \it Hilbert stable 
compactification \rm by Teixidor~\cite{Teix1}, of the universal
moduli space of semistable vector bundles of rank $r$ and degree $d$
over ${\frak M}_g$, the moduli space of smooth curves of genus $g$, 
by letting $C$ vary
and degenerate in ${\frak G}$. 
In contrast to the \it slope stable compactification \rm
of Pandharipande \cite{P} which involves torsion free sheaves on singular
curves, this compactification would take place entirely in the realm of
vector bundles. Its potential usefulness is illustrated by the paper
\cite{G2} where Hilbert stable vector bundles on a nodal curve
are used to prove a conjecture of Newstead and Ramanan on the moduli space
of stable rank two bundles over a smooth
curve. In order to make such a theory work, the objects one
starts with, namely Hilbert and Mumford stable vector bundles on
smooth curves, have to be same. Thus, one must show (1) that every stable
vector
bundle of rank $r$ over a smooth curve $C$ of sufficiently high degree $d$
gives rise to an embedding of $C$ into ${\frak G}$ and (2) that
for $C\hookrightarrow {\frak G}$ Hilbert and Mumford stability for
the bundle $E_C$ coincide. The first point follows from a recent theorem
of Butler~\cite{But} (see~\ref{But} below), and (2) has been established
in the rank two case by Gieseker and Morrison \cite{G1}.
It is the aim of the present note to settle the general case, i.e., prove
\begin{MThm}
Fix $g$ and $r$, then there is a constant $d_0$ such that for every
$d\ge d_0$ and every  complex vector space $W$ of dimension
$p=d+r(1-g)$ there exists a constant $m_0=m_0(d,g,r)$
such that for all $m\ge m_0$ the following holds true: 
Let $C\hookrightarrow G(W,r)$ be a smooth curve of genus $g$
and $W\otimes \O_C\lra E_C$ the restriction of the universal quotient
to $C$. Assume $W\lra H^0(E_C)$ is an isomorphism and $\deg(E_C)=d$. Then
$C$ is $m$-Hilbert (semi/poly)stable, if and only if $E_C$ is a
(semi/poly)stable vector bundle.
\end{MThm}
Note that both the condition of Mumford and Hilbert stability can
be formulated as stability requirements on the quotient $W\otimes\O_C\lra
E_C$. Therefore, it is a natural idea to look at the $\SL(W)$-action  on (some
open part of) the quot scheme of quotients of $W\otimes \O_C$.
As it turns out both stability conditions give rise to the same
linearized line bundle on this open part of the 
quot scheme. If the parameter space were projective, this would settle
the problem. Since this is not the case, we have to see how the
curve $C$ with $E_C$ semistable and Hilbert semistable
might degenerate in the set of Hilbert semistable points.
In turns out that the degeneration is roughly $C$ with some rational
components attached, a case which can be excluded by an adaptation of 
an argument
from~\cite{Teix1}. In other words, the locus of smooth curves $C^\p$ which
are isomorphic to $C$ such
that $E_{C^\p}$ is semistable is closed in the locus of Hilbert semistable
points. This is now as good as the projectivity of the parameter space
and one can conclude by standard methods in Geometric Invariant Theory.
Our proof therefore avoids completely any non-trivial computation. 
\section{Preliminaries}
\subsection{Review of some aspects of the theory of semistable vector
bundles}
A vector bundle over a smooth 
curve $C$ is called \it (semi)stable\rm, if it satisfies
$\mu(F)(\le)\mu(E)$ for all non-trivial proper subbundles $F\subset E$,
and \it polystable\rm, if $E$ is isomorphic to a direct sum of stable
bundles all of which have the same slope.
\par
The following is a recent generalization to semistable vector bundles
of a result of Mumford 
on line bundles.
\begin{Thm}[Butler~\cite{But}]
\label{But}
Let $E$ and $E^\p$ be semistable vector bundles on the smooth curve $C$ 
of genus $g$.
Assume $\mu(E)>2g$ and $\mu(E^\p)\ge 2g$. Then the homomorphism
$$
H^0(E)\otimes H^0(E^\p)\q\lra \q H^0(E\otimes E^\p)
$$
is surjective.
\end{Thm}
From this, one infers (see \cite{Teix3})
\begin{Cor}
\label{CorBut}
Let $E$ be a semistable vector bundle of rank $r$ on the smooth
curve $C$ of genus $g$ with $\mu(E)> 2g$. Then, the
homomorphism
$$
\bigwedge^r H^0(E)\q\lra\q H^0\bigl(\bigwedge^r E\bigr)
$$
is surjective.
\end{Cor}
Note that under the assumptions of~\ref{CorBut}, $H^1(E)=0$.
So, Corollary~\ref{CorBut} shows that the quotient
$H^0(E)\otimes\O_C\lra E$ defines an embedding of $C$ into
$G(W,r)$ where $W$ is a complex vector space of dimension $\deg(E)+r(1-g)$.
\begin{Prop}
\label{Sect}
Fix $g$ and $r$. Then there is a constant $d_1>2g$, such that for
every curve $C$ of genus $g$ and every vector bundle $E$ of rank $r$
and degree $d\ge d_1$ the following conditions are equivalent
\begin{enumerate}
\item $E$ is a (semi)stable vector bundle.
\item $h^0(F)/\rk F\ (\le)\ \chi(E)/r$ 
      for all non-trivial proper subbundles $F$ of $E$.
\item $\chi(E)/r\ (\le)\ h^0(Q)/\rk Q$ for all non-trivial proper quotient bundles
      $Q$ of $E$.
\end{enumerate}
\end{Prop}
\begin{pf} This is standard. See \cite{LP} or \cite{HL}. From the proof
one can easily determine an explicit value for $d_1$.
\end{pf}
\subsection{Properties of semistable points}
\label{InvTh}
Let $X$ be a quasi projective scheme on which the reductive group $G$ acts.
Suppose this action comes with a linearization in an ample line bundle
$A$. Then, the open sets $X^{ss}$ and $X^s$ of semistable and stable
points are defined. Furthermore, a semistable point $x$ is called
\it polystable\rm, if its orbit is closed in $X^{ss}$. The set
of polystable points will be denoted by $X^{ps}$.
Now, assume that $X$ is projective.
For any point $x\in X$ and any one parameter subgroup $\la$ of $G$,
we define $\mu_A(x,\la)$ as minus the weight of the $\C^*$-action
induced by $\la$ on the fibre of $A$ over the point 
$\lim_{x\rightarrow 0} \la(z)\cdot x$. The Hilbert-Mumford criterion
then says that a point $x$ is (semi)stable if and only if
$\mu_A(x,\la)(\ge)0$  holds for every one parameter subgroup $\la$ of $G$.
Moreover, $x$ is polystable if and only if it is semistable and
a fix point for every $\C^*$-action coming from a one parameter subgroup
$\la$ with $\mu_A(x,\la)=0$.
\par
Next, suppose we are given two representations $\rho_1\colon G\lra \SL(V_1)$
and $\rho_2\colon G\lra \SL(V_2)$ of the reductive group $G$ on the 
finite dimensional $\C$-vector spaces $V_1$ and $V_2$.
This yields an action of $G$ on $\P(V_1)\times \P(V_2)$ together with natural
linearizations in $\O(t_1,t_2)$ for all $t_1, t_2> 0$.
The corresponding set of (semi/poly)stable points depends only on the
parameter $\theta:= t_1/t_2\in (0,\infty)$ and will be denoted by
$Q_\theta^{(s/p)s}$. We also define $Q_0^{(s/p)s}$ and $Q_\infty^{(s)s}$ as the
preimage of the (semi/poly)stable points under the projection onto 
$\P(V_1)$ and $\P(V_2)$, respectively. Then, the following properties
are well known and easy to see (\cite{Th}, \cite{Re}): 
There exists a finite number of
critical values $\theta_1,...,\theta_s\in (0,\infty)$
such that,
settting $\theta_0=0$ and $\theta_{s+1}=\infty$,
for $i=1,...,s+1$ and given $\theta$, $\theta^\p$ in $(\theta_{i-1},\theta_i)$
\begin{eqnarray}
\label{eins}
Q_\theta^{(s/p)s} & =&Q_{\theta^\p}^{(s/p)s}\\
\label{zwei}
Q_\theta^{ss} &\subset & Q_{\theta_{i-1,i}}^{ss}\\
\label{drei}
Q_\theta^{s} &\supset & Q_{\theta_{i-1,i}}^s.
\end{eqnarray}
Now, let $X$ be a $G$-invariant closed subscheme of $\P(V_1)\times\P(V_2)$,
and set $X^{(s/p)s}_{\theta}:= Q_\theta^{(s/p)s}\cap X$, 
$\theta\in [0,\infty]$.
\begin{Lem}
\label{SemStab}
Suppose that there is an $n>0$, such that for every point $x\in X_0^{ss}$
and every one parameter subgroup $\la$ of $G$
$$
n\cdot \mu_{\O_{\P(V_1)}(1)}(\la, \pi_1(x))
\q\ge\q\mu_{\O_{\P(V_2)}(1)}(\la, \pi_2(x)).
$$
Then, for $\theta\in (0,\theta_1)$, also
$$
X_0^{(s)s}\q =\q X_\theta^{(s)s}.
$$
\end{Lem}
\begin{pf}
The stated condition clearly implies $X_\theta^s\subset X_0^s$, and thus,
by (\ref{drei}), $X_\theta^s= X_0^s$.
\par
There is a surjective morphism $X_\theta^{ss}\catqot G\stackrel{\phi}{\lra}
X_0^{ss}\catqot G$ which is by our assumption an isomorphism over
$X_0^{s}=X_\theta^s$. Thus, 
$X_\theta^{ss}\catqot G\setminus X_\theta^s\catqot G$
maps onto $X_0^{ss}\catqot G\setminus X_0^s\catqot G$ 
which means that for every
point $x\in X_0^{ss}\setminus X_0^s$ there exists a point $x^\p\in
X_\theta^{ss}\setminus X_\theta^s$ with $\phi([x^\p])=[x]$.
Choose $x^\p\in X_{\theta}^{ps}$. We claim that $x^\p$ also
lies in $X_0^{ps}$.
Indeed, let $\la$ be a one parameter subgroup with
$\mu_{\O_{\P(V_1)}(1)}(\la, \pi_1(x^\p))=0$. By the assumption and the
fact that $x^\p\in X_{\theta}^{ss}$, we must also have 
$\mu_{\O_{\P(V_2)}(1)}(\la, \pi_2(x^\p))=0$, and hence
$\mu_{\O(t_1,t_2)}(\la, x^\p)=0$ for all $t_1/t_2\in (0,\theta_1)$.
Since $x^\p$ is a fixed point for the corresponding $\C^*$-action, 
our claim is settled.
\par
Thus we have shown that for every $x\in X_0^{ss}$, the unique closed orbit
in $\overline{G\cdot x}$ is contained in $X_{\theta}^{ss}$,
whence also $G\cdot x\subset X_\theta^{ss}$ which is what we claimed.
\end{pf}
This argumentation also yields
\begin{Cor}
\label{SemStab2}
If, for $\theta\in (\theta_{i-1},\theta_i)$, one has
$X_{\theta}^{ps}\subset X_{\theta_{i-1}}^{ps}$, or
$X_{\theta}^{ps}\subset X_{\theta_{i}}^{ps}$, then
$X^{(s)s}_{\theta}=X_{\theta_{i-1}}^{(s)s}$, or
$X^{(s)s}_{\theta}=X_{\theta_{i}}^{(s)s}$, respectively.
\end{Cor}
\subsection{Some lemmas about Hilbert semistable curves}
In the rest of this paper, we will freely make use of the fact that,
if $C$ is a curve without embedded components, then the
restriction map $\E\lra \bigoplus_{i=1}^s\E_{C_i}$ is
injective for every locally free (or more generally depth 1) sheaf
$\E$ on $C$, where the $C_i$, $i=1,...,s$, are the components of $C$.
\par
Based on ideas of the papers \cite{G2} and \cite{Teix1},
we will now draw some consequences from the Hilbert
semistability of curves.
For this, fix $d$, $g$, and $r$ as before, and let ${\frak H}_{d,g}$
be the Hilbert scheme of all closed
subschemes of ${\frak G}$ with
Hilbert polynomial $P(m)=md+1-g$.
The notation $C\in {\frak H}_{d,g}$ means that $C$ is a closed
subscheme of ${\frak G}$ with Hilbert polynomial $P(m)$.
For any such $C$, the objects $E_C$ and $\psi^m_C$ are defined as
in the introduction.
First, since ${\frak H}_{d,g}$ is projective, and $\psi_C^1(\bigwedge^r W)
\subset H^0(\bigwedge^r E_C)$ is a very ample linear system, we can
find an $m_0^\p$, such that the map $\psi_C^m$ is surjective
for all $m\ge m_0^\p$ and for all $C\in {\frak H}_{d,g}$.
Hence, for $m\ge m_0^\p$ and $C\in {\frak H}_{d,g}$, 
the homomorphism $\phi_C^m$ is also defined,
and we may investigate the concept of $m$-Hilbert semistability for $C$. 
We set $V_2^m:=\bigwedge^{mP(m)} (S^m\bigwedge^r W)$.
\begin{Lem}
\label{Destab2}
There is an $m_0^{\p\p}\ge m_0^\p$, such that for every $m\ge m^{\p\p}_0$ and
every $C\in {\frak H}_{d,g}$ the following holds:
A subspace 
$
W_0\subset\ker\bigl( W\rightarrow H^0(E_{C|C_{\rm red}})\bigr)
$
gives rise to a one parameter subgroup $\la$ of $\SL(W)$
with 
$
\mu_{\O_{\P(V_2^m)}(1)}(\la, [\phi_C^m]) < 0.
$
Here, $C_{\rm red}$ stands for the reduced subscheme of $C$.
\end{Lem}
\begin{pf}
Choose a basis $v_1,...,v_{i_0}$ for $W_0$, complete it to a basis
$v_1,...,v_p$ of $W$, and define $\la$ w.r.t.\ this basis
by the weight vector $(i_0-p,...,i_0-p, i_0,...,i_0)$
where $i_0-p$ occurs $i_0$-times.
Set $W_1:=\langle\, v_{i_0+1},...,v_p\,\rangle$.
We obtain a splitting 
$\bigwedge^r W=\bigwedge^r W_1\oplus \widetilde{\Lambda}$.
The image of $\widetilde{\Lambda}$ in $H^0(\bigwedge^r E_C)$
lies in the kernel of the reduction
$H^0(\bigwedge^r E_C)\lra H^0(\bigwedge^r E_{C|C_{\rm red}})$,
in particular, the image
of $S^{m}\widetilde{\Lambda}$ in $H^0((\bigwedge^r E_C)^{\otimes m})$
is zero for all $m$ greater than some constant $\widetilde{m}$.
For those $m$, the minimum weight of an
eigenvector in $S^m \bigwedge^r W$
with non-zero image in $H^0((\bigwedge^r E_C)^{\otimes m})$
is $\widetilde{m}(i_0-p)+ (m-\widetilde{m})i_0
= m i_0 - \widetilde{m}p\ge m- \widetilde{m}p$, i.e.,
for $m>\widetilde{m}p$ we will definitely have
$\mu_{\O_{\P(V_2^m)}(1)}(\la, [\psi_C^m]) < 0$.
By the projectivity of ${\frak H}_{d,g}$, we can choose
$m_0^{\p\p}$ such that $m_0^{\p\p}\ge \widetilde{m}+1$ for every curve
$C\in {\frak H}_{d,g}$.
\end{pf}
The rest of this section will be devoted to prove a technical key result.
A curve $\widehat{C}$ will be called a \it tree-like curve\rm, if it satisfies
the following conditions
\begin{itemize}
\item $\widehat{C}$ is reduced, every reducible component is smooth, 
      meets at most
      two other components, and all intersections are ordinary double
      points.
\item The graph $\Gamma_{\widehat{C}}$ is a tree. Here, 
      $\Gamma_{\widehat{C}}$ is the graph
      with vertices $\{\, C_0,...,C_s\,\}$, the irreducible components
      of $\widehat{C}$, 
      and $C_i$ and $C_j$ are connected by an edge if and only if
      they meet.
\end{itemize}
We will call a vertex $C_i$ an \it end\rm , if there is only one edge
at $C_i$.
We will assume from now on that all irreducible components
of $\widehat{C}$ except $C_0$ are rational and that 
the genus of $\widehat{C}$ is $g$.
Suppose we are given a quotient $W\otimes\O_{\widehat{C}}\lra E$ 
where $E$ is a vector bundle of rank $r$ and degree $d$
and $\dim W=d+r(1-g)$. We label the vertex $C_i$
by $d_i:=\deg E_{|C_i}$, $i=0,...,s$, and set $d^\p:= d-d_0$.
Observe that, for $i\ge 1$, 
$E_{|C_i}\cong\O_{\P_1}(a_1)\oplus\cdots\oplus\O_{\P_1}(a_r)$
with $a_1\ge \cdots \ge a_r\ge 0$ and $\sum a_j = d_i$.
Suppose the induced homomorphism $W\lra H^0(E)$ is injective.
Let $H^0(E)\subset \bigoplus_{i=0}^m H^0(E_{|C_i})$ be the canonical
injection. Let $C_i$, $i\ge 1$, be an end (this exists),
and set $W_i:=\ker\bigl(W\rightarrow \bigoplus_{j\neq i} H^0(E_{|C_j}\bigr)$,
i.e., $W_i=W\cap H^0(E_{|C_i})\subset H^0(E_{|C_i}(-c_i))$, 
$c_i$ the point of intersection
of $C_i$ with the rest of the curve. Then, $\dim W_i\le d_i$,
i.e., $\dim(W/W_i)\ge d-d_i$. By removing
$C_i$ we obtain a new tree like curve $\widehat{C}^\p$ 
whose graph $\Gamma_{\widehat{C}^\p}$
is $\Gamma_{\widehat{C}}$ with the vertex $C_i$ and the edge at $C_i$ removed.
We can therefore iterate this procedure.
Set $W^\p:=
\ker\bigl(W\rightarrow \bigoplus_{i=1}^m H^0(E_{|C_i})\bigr)$.
If $\dim(W/W^\p)
= d-d^\p+ r(1-g)$, then we must have had equality at each step, whence
$\widetilde{W}:=\ker(W\rightarrow H^0(E_{|C_0}))$ identifies with 
$H^0(E_{\widetilde{C}}(-p_1-...-p_t))$ where $\widetilde{C}$ is the
closure of $\widehat{C}\setminus C_0$ in $\widehat{C}$, 
and $p_1,...,p_t$ are the points
of intersection of $C_0$ and $\widetilde{C}$.
\par
Next, consider the induced morphism $f^\p\colon \widehat{C}\lra G(W,r)$.
This morphism contracts all curves $C_i$ with $d_i=0$, in particular,
all ends labelled by $0$. For this reason, we can assume that
no such ends are present. The rational curves of positive degree
are embedded by $f^\p$, so that we can fix an end $C_{i_0}$
which is embedded by $f^\p$. The main result we will
need later is 
\begin{Prop}
\label{Teix}
There is an $m_0\ge m_0^{\p\p}$, such that, for every $m\ge m_0$,
the following conclusion is valid:
In the above situation, assume there is a curve $C^\p\in {\frak H}_{d,g}$
such that 
\begin{enumerate}
\item $\widehat{C}$ maps onto $C^\p_{\rm red}$, the reduction of $C^\p$.
\item $C^\p$ is generically reduced along the image of $C_{i_0}$.
\item The induced map $W\rightarrow H^0(E_{C^\p|C^\p_{\rm red}})$
      is injective.
\end{enumerate}
Let $C^\p_0$ be the component $f^\p(C_{i_0})$ of $C_{\rm red}^\p$,
and $R$ the union of the remaining components of $C^\p_{\rm red}$. 
Then we find a subspace $W_0\subset\ker\bigl(W\rightarrow 
H^0(E_{C^\p|R})\bigr)$,
such that, for a one parameter subgroup $\la$ of $\SL(W)$ associated to this
subspace as in the proof of Lemma~\ref{Destab2},
one gets $\mu_{\O_{\P(V_2^m)}(1)}
(\la, [\phi_{C^\p}^m]) < 0$.
\end{Prop}
\begin{pf}
There are canonical injective maps $\O_{C^\p_0}\subset f^\p_*\O_{C_{i_0}}$,
and $\O_R\subset f^\p_*\O_{\cup_{i\neq i_0}C_i}$.
For this reason and because of the third assumption,
the
maps 
$$
\begin{array}{ccl}
W &\lra & H^0(E_{C^\p|C^\p_{\rm red}})\\
  &\lra & H^0(E_{C^\p|C_0^\p})\oplus H^0(E_{C^\p|R})\\
  &\lra & H^0(C_{i_0},E_{|C_{i_0}})\oplus H^0(\bigcup_{i\neq i_0} C_i, 
                                              E_{|\cup_{i\neq i_0}C_i})
\end{array}
$$
are injective, whence $\ker(W\rightarrow H^0(E_{C^\p|R}))$ naturally
identifies
with $H^0(E_{|C_{i_0}}(-c_{i_0}))$.
Recall that $E_{|C_{i_0}}\cong \O_{\P_1}(a_1)\oplus\cdots\oplus\O_{\P_1}(a_r)$
with $a_1\ge \cdots \ge a_r\ge 0$ and $\sum a_j=d_{i_0}>0$,
whence $a_1\ge 1$.
We take $W_0=H^0(\O_{\P_1}(a_1-1))$ under these identifications.
\par
Let $\widetilde{C}_0^\p$ be the scheme theoretic closure
of the open subset $C^\p\setminus R$ in $C^\p$.
Define $\tau:=\dim\ker(H^0(\O_{\widetilde{C}_0^\p})
\rightarrow H^0(\O_{C_0^\p}))$.
Let $\L$ be an invertible sheaf on $C^\p$ and $\L^\p\subset \L$ a
subsheaf of $\L$ with support in $\widetilde{C}_0^\p$,
then
\begin{equation}
\label{100}
H^0(\L^\p_{|C_0^\p})\q \ge \q H^0(\L^\p)-\tau.
\end{equation}
\par
Now, we can apply the arguments used by Teixidor in \cite{Teix1},
Proof of~2.4.
Let $v_1,...,v_{j_0}$ be a basis for $W_0$, complete it
to a basis $v_1,...,v_p$ of $W$, and let $\la$ 
be given w.r.t.\ basis by $(\, j_0-p,...,j_0-p, j_0,...,j_0\,)$.
We also define $W_1:=\langle\, v_{j_0+1},...,v_p\,\rangle$.
The statement $\mu_{\O_{\P(V_2^m)}(1)}
(\la, [\phi_{C^\p}^m]) < 0$ can be translated into the statement
(cf.~\cite{G2}, \cite{Teix1})
\begin{equation}
\label{101}
{rmP(m)\over p} (p-a_1)\q <\q -\mu_{\O_{\P(V_2^m)}(1)}(\la^\p, 
[\phi_{C^\p}^m]).
\end{equation}
Here, $\la^\p$ is the one parameter subgroup of $\GL(W)$ given w.r.t.\
the fixed basis by the weight vector $(\, 0,...,0,1,...,1\,)$, $0$
appearing $j_0$-times.
Moreover,
\begin{eqnarray*}
-\mu_{\O_{\P(V_2^m)}(1)}(\la^\p, 
[\phi_{C^\p}^m]) & = & rmP(m)-\sum_{k=0}^{rm-1} b_k\\
                  & \ge & rmP(m)-\sum_{k=m(r-1)}^{rm -1} \tilde{b}_k
                          -rm\tau.
\end{eqnarray*}
Here, $b_k$ and $\tilde{b}_k$ are the dimensions of the subspaces
of $H^0((\bigwedge^r E_{C^\p})^{\otimes m})$ and
$H^0((\bigwedge^r E_{C^\p|C_0^\p})^{\otimes m})$ generated
by the eigenspace of weight $k$ in $S^m(\bigwedge^r W)$.
Note that only the space 
$S^m(W_0\otimes\bigwedge^{r-1} W_1\allowbreak \oplus \bigwedge^r W_1)$
yields non-zero sections in $H^0((\bigwedge^r E_{C^\p|C_0^\p})^{\otimes m})$,
so that the asserted inequality follows from (\ref{100}).
Next, by definition, for $m(r-1)\le k<mr$, the
image of the eigenspace of weight $k$
lies in $H^0(\P_1,\O_{\P_1}(md_{i_0}-(mn-k)))$, i.e.,
$\tilde{b}_k\le md_{i_0}+k-mn+1$.
The left hand side of (\ref{101}) is $m^2\cdot rd(1-a_1/p)+l_1(m)$,
$l_1(m)$ a linear polynomial, and the right hand side is bounded
from below by $m^2(rd-d_{i_0}+1/2)+l_2(m)-rm\tau$, $l_2(m)$ also a linear
polynomial.
Negating (\ref{101}) for large $m$ yields
$$
a_1\q\le\q {p\over rd}\left(d_{i_0}-\small{1\over 2}\right)\q 
<\q {d_{i_0}\over r},
$$
a contradiction.
\par
Now, the polynomials $l_1(m)$ and $l_2(m)$ depend only on $d,g,r$
and $d_{i_0}$ which leaves only finitely many possibilities
after fixing $d$, $g$, and $r$, because $0<d_{i_0}\le d$.
Moreover, $\tau$ is bounded by $h^0(\O_{C^\p})-1$, so it can take only finitely
many values for $C^\p$ varying in ${\frak H}_{d,g}$.
This means that we can indeed find $m_0$ as asserted.
\end{pf}
\section{Proof of the Theorem}
Choose $d>d_1$ according to Proposition~\ref{Sect}, 
fix a complex vector space $W$ 
of dimension $d+r(1-g)$, and let $C\hookrightarrow {\frak G}=G(W,r)$
be a smooth curve of genus $g$.
This provides us, on $C$, with a quotient $W\otimes\O_C\lra E_C$.
Write $L=L_C$ for the line bundle $\det E_C$.
Let ${\frak Q}_0$ be the quasi projective quot scheme parametrizing all
quotients $q\colon W\otimes \O_C\lra E$, such that
\begin{itemize}
\item $E$ is a vector bundle on $C$ of rank $r$
with determinant $L$ 
\item $H^0(q)$ is an isomorphism
\item $\bigwedge^rW\lra H^0(L)$ is surjective.
\end{itemize} 
\subsection{Review of Gieseker's construction of the moduli space
of stable bundles}
On ${\frak Q}_0\times C$, there is the universal quotient
$W\otimes \O_{{\frak Q}_0\times C}\lra {\frak E}_{{\frak Q}_0}$ which
provides us with $\bigwedge^r W\otimes \O_{{\frak Q}_0\times C}\lra
\bigwedge^r {\frak E}_{{\frak Q}_0}$. Note that
$\bigwedge^r {\frak E}_{{\frak Q}_0}\cong \pi_C^*L\otimes \pi_{{\frak Q}_0}^*
{\cal A}$ for some $\SL(W)$-linearized line bundle ${\cal A}$ on ${\frak Q}_0$,
so that projecting the latter homomorphism to ${\frak Q}_0$ yields
$$
\bigwedge^r W\otimes\O_{{\frak Q}_0}\lra H^0(L)\otimes {\cal A}.
$$
This homomorphism induces an injective and $\SL(W)$-equivariant
morphism $\iota\colon {\frak Q}_0\lra \P(V_1)$
with $V_1:=\Hom(\bigwedge^r W, H^0(L))^\vee$. Using Corollary~\ref{CorBut}
and Proposition~\ref{Sect},
it follows that the preimage under $\iota$ of the (semi/poly)stable
points is exactly the set of quotients $q\colon W\otimes \O_C\lra E$
for which $E$ is a (semi/poly)stable vector bundle.
Write ${\frak Q}^{(s/p)s}$ for the respective sets. The induced map
$\iota\colon {\frak Q}^{ss}\lra \P(V_1)^{ss}$ is proper, from which one infers
that ${\cal M}_{L/r}:={\frak Q}^{ss}\catqot \SL(W)$ exists. 
\subsection{Proof of the theorem}
Set $V_2^m:=\bigwedge^{P(m)} (S^m\bigwedge^rW)$, so that, for every $m\ge 1$, 
we have a natural morphism
$$
j_m\colon {\frak Q}_0 \lra \P(V_1)\times \P(V_2^m).
$$
\begin{Rem}
We remark in passing that the pullback of $\O(1)$
under the morphism ${\frak Q}_0\lra \P(V_2^m)$ is just
${\cal A}^{\otimes mP(m)}$, i.e., the morphisms 
from ${\frak Q}_0$ to $\P(V_1)$ and $\P(V_2^m)$ both give rise to the
same $\SL(W)$-linearized line bundle on ${\frak Q}_0$.
\end{Rem}
Let $X^m$ be the closure of $j_m({\frak Q}_0)$.
We will now use the notation of Section~\ref{InvTh}.
Note that for every point 
$x=([x_1], [x_2])=j_m([q\colon W\otimes \O_C\lra E])$ and
every one parameter subgroup $\la$ of $\SL(W)$, we have 
\begin{equation}
\label{sechs}
m P(m) \cdot \mu_{\O_{\P(V_1)}(1)}(\la, [x_1])\q \ge\q 
\mu_{\O_{\P(V_2^m)}(1)}(\la, [x_2]),
\end{equation}
so that in view of Prop.~\ref{Sect},
one immediately infers
\begin{Cor} If the curve $C$ is $m$-Hilbert (semi)stable, then the
vector bundle $E_C$ is (semi)stable.
\end{Cor}
\begin{Rem}
Note that this conclusion holds for every $m\ge 1$.
\end{Rem}
By Lemma~\ref{SemStab},
$
j_m({\frak Q}^{ss}) = (X^m)^{ss}_\theta\subset  (X^m)_{\theta_1}^{ss}.
$
Suppose now we could prove the following
\begin{Prop}
\label{Crucial}
Let $m_0$ be as in Prop.~\ref{Teix}. 
Then for all $m\ge m_0$, the following holds true:
Let $[q\colon W\otimes \O_C\lra E]\in {\frak Q}^{ss}$, such that
$E$ is a polystable vector bundle. Suppose $j_m([q])\in (X^m)^{ss}_{\theta}$
for some $\theta\in (0,\infty)$.
Then, also
$$
j_m([q])\q \in\q (X^m)^{ps}_\theta.
$$
\end{Prop}
In this case, by Corollary~\ref{SemStab2},
$
j_m({\frak Q}^{(s)s}) =  (X^m)^{(s)s}_{\theta_1}.
$
Using (\ref{sechs}) and Lemma~\ref{SemStab} again,
we also get
$
(X^m)^{(s)s}_{\theta}  = (X^m)_{\theta_1}^{(s)s}$
for all
$\theta\in (\theta_1,\theta_2)$.
Now, iterating this argumentation, yields the conclusion
$$
(X^m)_\infty^{(s)s}\q =\q j_m({\frak Q}^{(s)s})
$$
which is just a reformulation of the assertion of the theorem.\qed
\subsection{Proof of Proposition~\ref{Crucial}}
Let $\la$ be a one parameter subgroup of $\SL(W)$, such that
$x_0:=\lim_{z\rightarrow 0} j_m([q])\cdot \la(z)$ exists
in $(X^m)^{ss}_\theta$, but such that $j_m([q])$ is not a fixed point
for the corresponding $\C^*$-action.
We must describe $x_0=([x_1],[x_2])$ more explicitly to derive
a contradiction. First, by assumption,
we have a morphism $\C^*\lra {\frak Q}^{ss}$.
This corresponds to a family $W\otimes \O_{\C^*\times C}\lra {\E}_{\C^*}$.
This family can be extended to a family of quotients
$W\otimes\O_{\C\times C}\lra \E_\C$ where $\E_\C$ is a $\C$-flat
family of coherent sheaves of rank $r$ with determinant $L$ on $C$.
Note that the flatness over $\C$ implies that $\E_\C$ is torsion free
as $\O_{\C\times C}$-module.
Set $E_\C:=\E_\C^{\vee\vee}$. This is a reflexive sheaf
on the smooth surface $\C^*\times C$, whence it is locally free and
thus flat over $\C$. This gives a 
family
$$
W\otimes \O_{\C\times C}\lra E_\C.
$$
\begin{Rem}
\label{remind}
Let us remind the reader of some features of this construction.
\begin{enumerate}
\item The kernel of the homomorphism 
$
\E_{\C|\{0\}\times C} \lra E_{\C|\{0\}\times C}
$
is exactly the torsion ${\cal T}$ of $\E_{\C|\{0\}\times C}$.
\item Since $W\otimes \O_C$ generically generates $ E_{\C|\{0\}\times C}$,
we see $\dim_\C({\cal T})\le\dim W-r=d-rg<d$, thus 
$\deg(\E_{\C|\{0\}\times C}/{\cal T})>0$
has positive degree, and since there is a surjection
$W\otimes \O_C\lra \E_{\C|\{0\}\times C}/{\cal T}$, the rational map
$C \dasharrow {\frak G}$ induced by $W\otimes C\lra E_{\C|\{0\}\times C}$
is not constant.
\item Set $\widetilde{W}:=W/\ker(W\rightarrow H^0(\E_{\C|\{0\}\times C}))$.
Then, $\dim\Im(W\rightarrow \allowbreak H^0(E_{\C|\{0\}\times C}))=
\dim(\widetilde{W})-d^\p$, by 1.
\end{enumerate}
\end{Rem}
From this discussion, we deduce that the homomorphism
$$
\bigwedge^r W\otimes \O_{\C\times C}\lra \bigwedge^r E_\C
$$
is surjective outside a finite set of points $p_1,...,p_t$
located on $\{0\}\times C$ where $t\le d^\p$.
In particular, there is a rational map
$$
h\colon \C\times C\dasharrow {\frak G}
$$
defined outside $\{\, p_1,...,p_t\,\}$.
By blowing up the points $p_1,...,p_t$ and possibly some infinitely near
ones (see \cite{Beau}, II.7), we arrive
at a smooth surface $\widehat{S}$ together
with a morphism
$$
\widehat{h}\colon \widehat{S}\lra {\frak G}.
$$
\begin{Rem}
\label{Destab1}
The map $\bigwedge^r W\otimes\O_{\{0\}\times C}\lra H^0(L)=
\allowbreak 
H^0(\det  E_{\C|\{0\}\times C})$ defines the point $[x_1]\in \P(V_1)$.
As in the proof of~\ref{Destab2},
every subspace $W_0\subset \ker(W\rightarrow H^0( E_{\C|\{0\}\times C}))$
yields a one parameter subgroup $\la$ of $\SL(W)$
with $\mu_{\O_{\P(V_1)}(1)}(\la, [x_1])<0$.
\end{Rem}
The composite morphism $\widehat{S}\lra \C$ is still flat, and an easy
inductive argument shows that the
fibre $\widehat{C}$ over 
$\{0\}$ is a tree-like curve with $C$ as its only non-rational
component.
\par
Next, observe that by Butler's results \ref{CorBut}, the
morphism
$\C^*\times C\lra\C^*\times {\frak G}
\hookrightarrow \C^*\times\P(\bigwedge^r W)$
is an embedding and consequently corresponds to
a morphism $\C^*\lra {\frak H}_{d,g}$. By extending this morphism
to a morphism $\C\lra {\frak H}_{d,g}$, we get another
surface $S^\p$ equipped with a flat morphism to $\C$.
Observe that the flatness over $\C$ together with the fact that
$S^\p\times_\C\C^*$ is integral implies that $S^\p$ is also integral
(\cite{Ha}, III, Prop.9.7). Moreover, by our construction,
there a morphsim $f\colon \widehat{S}\lra S^\p$ which factorizes over
$\widehat{S}\lra \widetilde{S}^\p$, $\widetilde{S}^\p$ 
the normalization of $S^\p$.
The latter morphism just being the contraction of some rational curves
with negative self-intersection, the morphism $f\colon \widehat{S}\lra S^\p$
is proper.
\par
Now, write $\C=\Spec\C[T]$, and denote by $T$ also the induced element
in the function field $K(\widehat{S})=K(S^\p)$. 
Let $C^\p$ be the fibre of $S^\p$ over $\{0\}$. We will have to compare
$\widehat{C}$ and $C^\p$. For this let $[\widehat{C}]$ and
$[C^\p]$ be the Weil divisor classes of those curves.
By definition $\widehat{C}$ and $C^\p$ are the Cartier divisors
$\mathop{\rm div}(T)$, taken on $\widehat{S}$ and $S^\p$, respectively.
Proposition~1.4 in \cite{Fu} thus shows that $f_*[\widehat{C}]=[C^\p]$
on the cycle level. The upshot of this discussion is that,
if we can show that every rational curve in $\widehat{C}$
which is not contracted is mapped injectively to ${\frak G}$,
the only component of $C^\p$ which is possibly not generically reduced is 
$f(C)$, the ultimate goal being to apply Proposition~\ref{Teix}.
\par
Anyway, at this stage we know that the curve $C^\p\in{\frak H}_{d,g}$
supplies $[x_2]$ in $\P(V_2^m)$. Therefore, we can look at some destabilizing
one parameter subgroups.
\begin{Lem}
\label{Destab3}
For $m\ge m_0$,
the homomorphism $W\lra H^0(E_{C^\p_{\rm red}})$ must be injective.
\end{Lem}
\begin{pf}
Observe that a subspace $W_0$ of $\ker(W\rightarrow  H^0(E_{C^\p_{\rm red}}))$
gives by Lemma~\ref{Destab2} and Remark~\ref{Destab1} rise
to a one parameter subgroup $\la$ with
both $\mu_{\O_{\P(V_1)}(1)}(\la, [x_1])<0$ and 
 $\mu_{\O_{\P(V_2^m)}(1)}(\la, [x_2])<0$, in contradiction to
the semistability of $x_0$.
\end{pf}
The induced morphism $f^\p\colon \widehat{C}\lra C^\p_{\rm red}$
is surjective, so that there are
injections $\O_{C^\p_{\rm red}}\subset f^\p_*\O_{\widehat{C}}$,
and $E_{C^\p|C^\p_{\rm red}}\subset E_{C^\p|C^\p_{\rm red}}\otimes
f^\p_*\O_{\widehat{C}}=\allowbreak
f^\p_* f^{\p*} E_{C^\p|C^\p_{\rm red}}$.
The composite $W\lra H^0(E_{C^\p|C^\p_{\rm red}})\subset H^0(f^{\p*} 
E_{C^\p|C^\p_{\rm red}})$
is thus injective by Lemma~\ref{Destab3}.
Therefore, in Remark~\ref{remind}, 3., the space $\widetilde{W}$
equals $H^0({\cal T})$ and, thus, has dimension $d^\p$. 
Now, one immediately checks that we 
are exactly in the position to apply Proposition~\ref{Teix}.
Since the subspace $W_0$ used to destabilize $[x_2]$ lies
in the kernel of $W\rightarrow H^0(E_{\C|\{0\}\times C})$,
we find again a one parameter subgroup $\la$ of $\SL(W)$ with
$\mu_{\O_{\P(V_1)}(1)}(\la, [x_1])<0$ and 
$\mu_{\O_{\P(V_2^m)}(1)}(\la, [x_2])\allowbreak <0$, contradicting the assumptions 
made on $x_0$.\qed

\end{document}